\documentclass[noinfoline]{imsart}
\usepackage{latexsym,epsfig,amssymb,amsmath,amsthm,color,url,bbm,pdfsync}
\usepackage{graphicx,float,subfigure,bm}
\usepackage{amsfonts}
\usepackage{graphics}
\usepackage{mathrsfs,hyperref}

\numberwithin{equation}{section}
\allowdisplaybreaks
\usepackage[round,comma]{natbib}

\newtheorem{lemma}{Lemma}[section]

\newtheorem{theorem}[lemma]{Theorem}
\newtheorem{proposition}[lemma]{Proposition}
\newtheorem{definition}[lemma]{Definition}
\newtheorem{corollary}[lemma]{Corollary}
\newtheorem{example}[lemma]{Example}
\newtheorem{exercise}[lemma]{Exercise}
\newtheorem{remark}[lemma]{Remark}

\newcommand{\cov}{{\rm cov}}
\newcommand{\var}{{\rm Var}}
\newcommand{\cid}{\stackrel{d}{\rightarrow}}
\newcommand{\cip}{\stackrel{p}{\rightarrow}}

\newcommand{\stas}{\stackrel{\rm a.s.}{\rightarrow}}
\newcommand{\ex}{{\rm e}\,}

\newcommand{\bbr}{{\mathbb R}}
\newcommand{\bbz}{{\mathbb Z}}

\newcommand{\bfX}{{\bf X}}

\newcommand{\E}{{\mathbb E}}
\renewcommand{\P}{{\mathbb P}}

\newcommand{\beao}{\begin{eqnarray*}}
\newcommand{\eeao}{\end{eqnarray*}\noindent}
\newcommand{\beam}{\begin{eqnarray}}
\newcommand{\eeam}{\end{eqnarray}\noindent}
\newcommand{\beqq}{\begin{equation}}
\newcommand{\eeqq}{\end{equation}\noindent}
\newcommand{\beqo}{\begin{equation*}}
\newcommand{\eeqo}{\end{equation*}\noindent}
\newcommand{\bth}{\begin{theorem}}
\newcommand{\ethe}{\end{theorem}}

\newcommand{\bre}{\begin{remark}\em }
\newcommand{\ere}{\end{remark}}

\newcommand{\ble}{\begin{lemma}}
\newcommand{\ele}{\end{lemma}}

\newcommand{\bde}{\begin{definition}}
\newcommand{\ede}{\end{definition}}
\newcommand{\bco}{\begin{corollary}}
\newcommand{\eco}{\end{corollary}}
\newcommand{\bpr}{\begin{proposition}}
\newcommand{\epr}{\end{proposition}}
\newcommand{\bpf}{\begin{proof}}
\newcommand{\epf}{\end{proof}}
\newcommand{\bexer}{\begin{exercise}}
\newcommand{\eexer}{\end{exercise}}
\newcommand{\bexam}{\begin{example}\rm }
\newcommand{\eexam}{\end{example}}
\newcommand\1{{\mathbf 1}}

\newcommand{\bbeta}{{\boldsymbol\beta}}

\newcommand{\btheta}{{\boldsymbol\theta}}

\newcommand{\bSig}{{\boldsymbol\Sigma}}
\newcommand{\bTheta}{{\boldsymbol\Theta}}

\newcommand{\bu}{{\mathbf{u}}}
\newcommand{\boldm}{{\mathbf{m}}}
\newcommand{\bB}{{\mathbf{B}}}
\newcommand{\bA}{{\mathbf{A}}}

\newcommand{\bl}{{\mathbf{L}}}
\newcommand{\bQ}{{\mathbf{Q}}}
\newcommand{\calF}{{\mathcal{F}}}

\newcommand\independent{\protect\mathpalette{\protect\independenT}{\perp}}
\def\independenT#1#2{\mathrel{\rlap{$#1#2$}\mkern2mu{#1#2}}}

\begin{document}



\begin{frontmatter}
\title{\Large Goodness-of-Fit Testing for Time Series Models via Distance Covariance}
\runtitle{Goodness-of-Fit Testing for Time Series Models via Distance Correlation}
\begin{aug}
\author{Phyllis Wan}\ead[label=e1]{wan@ese.eur.nl}
\address{Econometric Institute, Erasmus University Rotterdam, \\P.O. Box 1738, 
3000DR Rotterdam, \\The Netherlands\\
\printead{e1}}
\author{Richard A. Davis}
\ead[label=e2]{rdavis@stat.columbia.edu}
\address{Department of Statistics, Columbia University, \\1255 Amsterdam Avenue, MC 4690, \\New York, NY 10027, USA\\
\printead{e2}}

\runauthor{Wan and Davis}


\end{aug}


\begin{abstract}
In many statistical modeling frameworks, goodness-of-fit tests are typically administered to the estimated residuals.  In the time series setting, whiteness of the residuals is assessed using the sample autocorrelation function.  For many time series models, especially those used for financial time series, the key assumption on the residuals is that they are in fact independent and not just uncorrelated.  In this paper, we apply the auto-distance covariance function  (ADCV) to evaluate the serial dependence of the estimated residuals.  Distance covariance can discriminate between dependence and independence of two random vectors.   The limit behavior of the test statistic based on the ADCV is derived for a general class of time series models.  One of the key aspects in this theory is adjusting for the dependence that arises due to parameter estimation.  This adjustment has essentially the same form regardless of the model specification.   We illustrate the results in simulated examples.

\end{abstract}

\end{frontmatter}
\maketitle
\bigskip


\section{Introduction}\label{sec:intro}

Let $\{X_j,j\in\bbz\}$ be a stationary time series of random variables with finite mean and variance. Given consecutive observations of this time series $X_1,\ldots,X_n$, we consider testing the plausibility that the data were generated from a parametric model.  We consider causal models of the form
\beqq \label{eq:model}
	X_j = g(Z_{-\infty:j};\bbeta),
\eeqq
where the $Z_j$'s are independent and identically distributed (iid) with mean zero and finite variance, $Z_{n_1:n_2}$ denotes the sequence $\{Z_j,n_1\le j \le n_2\}$, and $\bbeta\in\bbr^d$ is the parameter vector.  
Assume further that the model \eqref{eq:model} has the invertible representation
\beqq \label{eq:invert}
	Z_j  = h(X_{-\infty:j};\bbeta).
\eeqq
The objective of this paper is to provide a validity check of the model \eqref{eq:model} by testing the estimated residuals for {\it independence}.

Given observations $X_{1:n}$ and $\hat\bbeta$, an estimator for $\bbeta$, the innovations $\{Z_j\}$ can be approximated by the residuals based on the infinite sequence $X_{-\infty:j }$, defined as
\beqq \label{eq:tildez}
	\tilde Z_j := Z_j(\hat\bbeta) = h(X_{-\infty:j };\hat\bbeta).
\eeqq
Since we do not observe $X_j$ for $j\le0$, we instead use the estimated residuals
\beqq \label{eq:res}
	\hat Z_j := h(Y_{-\infty:j};\hat\bbeta),\quad j=1,\ldots,n,
\eeqq
where $\{Y_j\}$ is the infinite sequence with $Y_j=X_j$, $j \ge 1$ and $Y_j=0$ for $j\le 0$.
If the time series $\{X_j\}$ is stationary and ergodic, the influence of $X_{-\infty:0}$ in \eqref{eq:tildez} becomes negligible for large $j$ and $\hat{Z}_j$ and $\tilde{Z}_j$ become close.

It is general practice to inspect $\{\hat Z_j\}$ for goodness-of-fit of the time series model.
If \eqref{eq:model} correctly describes the generating mechanism of $\{X_j\}$, one would expect $\{\hat Z_j\}$ to behave similarly as $\{Z_j\}$.  However, the sequence $\{\hat{Z}_j,1\le j \le n\}$ is not iid since they are functions of $\hat\bbeta$, hence certain properties of $\{\hat Z_j\}$ can differ from that of $\{Z_j\}$, which in turn may impact sample statistics such as the sample autocorrelation of the residuals.  This has been noted for specific time series models in the literature.  For example, for the ARMA model, corrections have been made for statistics based on the residuals, see Section 9.4 of \cite{brockwell:davis:1991}.  For heteroscedastic GARCH models, the moment sum process of the residuals is notedly different from that of iid innovations, see \cite{kulperger:yu:2005}.  Though $\{\hat{Z}_j\}$ should be nearly independent under the true model assumption,  the discrepancy between $\{\hat{Z}_j\}$ and $\{{Z}_j\}$ should be taken into account when designing a goodness-of-fit test.  

In this paper, we characterize the serial dependence of the residuals using distance covariance.  Distance covariance is a useful dependence measure with the ability to detect both linear and nonlinear dependence.  It is zero if and only if independence occurs.  We study the auto-distance covariance function (ADCV) of the residuals and derive its limit when the model is correctly specified.  We show that the limiting distribution of the ADCV of $\{\hat Z_j\}$ differs from that of its iid counterpart $\{Z_j\}$ and quantify the difference.  This is an extension of Section 4 of \cite{davis:matsui:mikosch:wan:2018} which considered this problem for AR processes.

The remainder of the paper is structured as follows.  An introduction to distance correlation and ADCV along with some historical remarks are given in Section~\ref{sec:dcor}.  In Section~\ref{sec:meta}, we provide the limit result for the ADCV of the residuals for a general class of time series models.  To implement the limiting results, we apply the parametric bootstrap, the methodology and thoeretical justification of which is given in Section~\ref{sec:boot}.  We then apply the result to ARMA and GARCH models in Sections~\ref{sec:arma} and \ref{sec:garch} and illustrate with simulation studies.  A simulated example where the data does not conform with the model is demonstrated in Section~\ref{sec:noncausal}.


\section{Distance covariance} \label{sec:dcor}

Let $X\in\bbr^p$ and $Y\in\bbr^q$ be two random vectors, potentially of different dimensions.  Let  $\varphi_{X,Y}(s,t),\varphi_X(s),\varphi_Y(t)$ denote the joint and marginal characteristic functions of $(X,Y)$. We know that
$$
	X \independent Y \quad \Longleftrightarrow \quad \varphi_{X,Y}(s,t) = \varphi_X(s)\,\varphi_Y(t).
$$  
The {\em distance covariance} between $X$ and $Y$ is defined as
\beqo \label{eq:dcov}
	T(X,Y;\mu) = \int_{\bbr^{p+q}} \big|\varphi_{X,Y}(s,t)-\varphi_X(s)\,\varphi_Y(t)\big|^2\,\mu(ds,dt)\,,\,\quad (s,t)\in\bbr^{p+q},
\eeqo
where $\mu$ is a suitable measure on $\bbr^{p+q}$.
In order to ensure that $T(X,Y;\mu)$ is well-defined, one of the following conditions is assumed to be satisfied \citep{davis:matsui:mikosch:wan:2018}:
\vspace{-.05in}
\begin{enumerate}
\item
	$\mu$ is a finite measure;
\item
	$\mu$ is an infinite measure such that
	$$
		\int_{\bbr^{p+q}} (1\wedge|s|^\alpha) (1\wedge|t|^\alpha) \mu(ds,dt) < \infty \quad \text{and} \quad \E[|XY|^\alpha+|X|^\alpha+|Y|^\alpha] < \infty, \quad \text{for some $\alpha \in (0,2]$.}
	$$
\end{enumerate}
If  $\mu$ has a positive Lebesgue density on $\bbr^{p+q}$, then $X$ and $Y$ are independent if and only if $T(X,Y;\mu)=0$. 

For a stationary series $\{X_j\}$, the {\em auto-distance covariance} (ADCV) is given by
$$
	T_h(X;\mu) := T(X_0,X_h;\mu) = \int_{\bbr^2} \big|\varphi_{X_0,X_h}(s,t)-\varphi_X(s)\,\varphi_X(t)\big|^2\,\mu(ds,dt)\,,\,\quad (s,t)\in\bbr^2.
$$
Given observations $\{X_j,1\le j \le n\}$, the ADCV can be estimated by its sample version
$$
	\hat T_h(X;\mu) := \int_{\bbr^2} \big|C_n^X(s,t)\big|^2\,\mu(ds,dt)\,,\,\quad (s,t)\in\bbr^2,
$$
where 
$$
	C_n^{X}(s,t) := \frac{1}{n} \sum_{j=1}^{n-h} \ex^{isX_j+itX_{j+h}}-\frac{1}{n} \sum_{j=1}^{n-h} \ex^{isX_j} \frac{1}{n}  \sum_{j=1}^{n-h} \ex^{itX_{j+h}}.
$$
If we assume that $\mu = \mu_1\times\mu_2$ and is symmetric about the origin, {then under the conditions where $T_h(X;\mu)$ exists,} $\hat T_h(X;\mu)$ is computable in an alternative expression similar to a $V$-statistic, see Section~2.2 of \cite{davis:matsui:mikosch:wan:2018} for details.
It can be shown that if the $X_j$'s are iid, the process $\sqrt{n}C_n^X(s,t)$ converges weakly,
\beqq \label{eq:gh}
	\sqrt{n}C_n^X \cid G_h \quad \text{on $\mathcal{C}(K)$},
\eeqq
for any compact set $K \subset \bbr^2$, and
$$
	 n\hat T_h(X;\mu) \cid \int|G_h|^2 \mu(ds,dt),
$$
where $G_h$ is a zero-mean Gaussian process 
with covariance structure
\beao
\Gamma((s,t),(s',t')) &=& \cov(G_h(s,t),G_h(s',t')) \nonumber \\
&=& \E\big[\big(\ex^{i\,\langle s, X_0\rangle} - \varphi_X(s)\big)\big(\ex^{i\,\langle t,X_{h}\rangle} -
\varphi_X(t)\big) \nonumber \\
&&\hspace{0.5cm} \times \big(\ex^{-i\,\langle s',X_0\rangle} -
\varphi_X(-s')\big)\big(\ex^{-i\,\langle t',X_{h}\rangle} - \varphi_X(-t')\big)\big]\,.
\eeao


The concept of distance covariance was first proposed by \cite{feuerverger:1993} in the bivariate case and later popularized by \cite{szekely:rizzo:bakirov:2007}.  The idea of ADCV was first introduced by \cite{zhou:2012}.  For distance covariance in the time series context, we refer to \cite{davis:matsui:mikosch:wan:2018} for theory in a general framework.  

Most literature on distance covariance focus on the specific weight measure $\mu(s,t)$ with density proportional to $|s|^{-p-1}|t|^{-q-1}$.  This distance covariance has the advantage of being scale and rotational invariant, but imposes moment constraints on the variables under consideration.  In our case, as will be shown in Section~\ref{sec:meta}, this measure may not work when applied to the residuals (see also Section 4 of \cite{davis:matsui:mikosch:wan:2018} for a counterexample).  To avoid this difficulty, we assume a finite measure for $\mu$.  In this case $\hat T_h(X;\mu)$ has the computable form
\beao
	\hat T_h(X;\mu) &=& \frac{1}{(n-h)^2}\sum_{i,j=1}^{n-h} \hat\mu(X_i-X_j,X_{i+h}-X_{j+h}) \\
	&& \quad+ \frac{1}{(n-h)^4}\sum_{i,j,k,l=1}^{n-h} \hat\mu(X_i-X_j,X_{k+h}-X_{l+h}) \\
	&&\qquad - 2\frac{1}{(n-h)^3}\sum_{i,j,k=1}^{n-h} \hat\mu(X_i-X_j,X_{i+h}-X_{k+h}),
\eeao
where $\hat\mu(x,y)=\int \exp(isx+ity)\mu(ds,dt)$ is the Fourier transform with respect to $\mu$.

It should be noted that the concept of distance covariance is closely related to the Hilbert-Schmidt Independence Criterion (HSIC), see \cite{gretton:bousquet:smola:scholkopf:2005}.  For example, the distance covariance with Gaussian measure coincides with the HSIC with a Gaussian kernel.  In recent work, \cite{wang:li:zhu:2018} use HSIC for testing the cross dependence between two time series.


\section{General result} \label{sec:meta}

Let $X_1,\ldots, X_n$ be observations from a stationary time series $\{X_j\}$ generated from \eqref{eq:model} with $\beta=\beta_0$.  Let $\hat{Z}_1,\ldots,\hat{Z}_n$ be the estimated residual calculated through \eqref{eq:res}.  In this section, we examine the ADCV of the residuals
\beqo\label{eq:adcv}
	\hat{T}_h(\hat Z;\mu) := \|C_n^{\hat Z}\|^2_\mu = \int |C_n^{\hat Z}|^2 \mu(ds,dt),
\eeqo
where
$$
	C_n^{\hat Z}(s,t) := \frac{1}{n} \sum_{j=1}^{n-h} \ex^{is\hat Z_j+it\hat Z_{j+h}}-\frac{1}{n} \sum_{j=1}^{n-h} \ex^{is\hat Z_j} \frac{1}{n}  \sum_{j=1}^{n-h} \ex^{it\hat Z_{j+h}}.
$$
To provide the limiting result for $\hat{T}_h(\hat Z;\mu)$, we require the following assumptions.

\begin{enumerate}
\item[(M1)]\label{cond:m1}
	Let $\calF_j$ be the $\sigma$-algebra generated by $\{X_k,k\le j\}$.  We assume that the parameter estimate $\hat\bbeta$ is of the form
	\beqq\label{eq:m11}
		\sqrt{n}(\hat\bbeta - \bbeta_0) = \frac{1}{\sqrt{n}} \sum_{j=1}^n \boldm(Z_{-\infty:j};\bbeta_0) + o_p(1),
	\eeqq
	where $\boldm$ is a vector-valued function of the infinite sequence $X_{-\infty:j}$ such that
	\beqq\label{eq:m12}
		\E [\boldm(Z_{-\infty:j};\bbeta_0)|\calF_{j-1}] = \mathbf0, \quad \E|\boldm(Z_{-\infty:0};\bbeta_0)|^2 < \infty.
	\eeqq
	This representation can be readily found in most likelihood-based estimators, for example, the Yule-Walker estimator for AR processes, quasi-MLE for GARCH processes, etc.  In these cases $\boldm$ can be taken as the likelihood score function. By the martingale central limit theorem, \eqref{eq:m11} and \eqref{eq:m12} imply that
	$$
		\sqrt{n}(\hat\bbeta - \bbeta_0) \cid \bQ,
	$$
for a random Gaussian vector $\bQ$.
\item[(M2)]\label{cond:m2}
	Assume that the function $h$ in the invertible representation \eqref{eq:invert} is continuously differentiable, and writing
	\beqq \label{eq:bigL}
		\bl_j(\bbeta) := \frac{\partial}{\partial\bbeta} h(X_{-\infty:j};\bbeta),
	\eeqq
	we assume
	\beqo\label{eq:3}
		\E\|\bl_0(\bbeta_0)\|^2 <\infty.
	\eeqo
\item[(M3)]\label{cond:m3}
	Assume that $\{\hat Z_j\}$, the estimated residuals based on the {\em finite} sequence of observations,  is close to $\{\tilde Z_j\}$,  the fitted residuals based on the {\em infinite} sequence, such that
	$$
		\frac{1}{\sqrt{n}} \sum_{j=1}^n |\hat{Z}_j - \tilde{Z}_j|^k = o_p(1), \quad k=1,2.
	$$
\end{enumerate}


\bth \label{thm:meta}
Let $X_1,\ldots,X_n$ be a sequence of observations generated from the causal and invertible time series model \eqref{eq:model} and \eqref{eq:invert} with $\bbeta=\bbeta_0$. Let $\hat\bbeta$ be an estimator of $\bbeta$ and let $\hat Z_1,\ldots,\hat Z_n$ be the estimated residuals calculated through \eqref{eq:res} satisfying conditions \hyperref[cond:m1]{(M1)}--\hyperref[cond:m3]{(M3)}.  Furthermore assume that the weight measure $\mu$ satisfies
\beqq\label{eq:7}
	\int_{\bbr^2}\big[(1\wedge |s|^2)\,(1\wedge |t|^2) +  (s^2+t^2)\,\1(|s| \wedge |t|>1)\big]\mu(ds,dt) < \infty.
\eeqq
Then 
\beqo \label{eq:adcv:limit}
	n \hat{T}_h(\hat Z;\mu) \cid \|G_h + \xi_h\|^2_\mu,
\eeqo
where $G_h$ is the limiting distribution for $n \hat{T}_h(Z;\mu)$, the ADCV based on the iid innovations $Z_1,\ldots,Z_n$, and the correction term $\xi_h$ is given by
\beqq \label{eq:xi}
	\xi_h(s,t) := it\bQ^T \E\left[ \left(e^{isZ_0}-\varphi_Z(s)\right) e^{itZ_h} \bl_h(\bbeta_0)\right],
\eeqq
with $\bQ$ being the limit distribution of $\sqrt{n}(\hat\bbeta-\bbeta_0)$ and $\bl_h$ as defined in \eqref{eq:bigL}.
\ethe
The proof of the theorem is provided in Appendix~\ref{app:meta}.  
%

\bre

{\em Distance correlation}, analogous to linear correlation, is the normalized version of distance covariance, defined as
$$
	R(X,Y;\mu) := \frac{T(X,Y;\mu)}{\sqrt{T(X,X;\mu)T(Y,Y;\mu)}} \in [0,1].
$$
The {\em auto-distance correlation function} (ADCF) of a stationary series $\{X_j\}$ at lag $h$ is given by
$$
	R_h(X;\mu) := R(X_0,X_h;\mu),
$$
and its sample version $\hat R_h(X;\mu)$ can defined similarly.  It can be shown that the ADCF for the residuals from an AR($p$) model has the limiting distribution \citep{davis:matsui:mikosch:wan:2018}:
\beqq \label{eq:adcf:limit}
	n \hat{R}_h(\hat Z;\mu) \cid \frac{\|G_h + \xi_h\|^2_\mu}{T_0(Z;\mu)},
\eeqq
and the result can be easily generalized to other models.
In the examples in Sections~\ref{sec:arma} and \ref{sec:garch}, we shall use ADCF in place of ADCV.

\ere


\section{Parametric bootstrap}\label{sec:boot}

The limit in \eqref{eq:adcf:limit} is not distribution-free and is generally intractable.  In order to use the result, we propose to approximate the limit through the parametric bootstrap described below.

Given observations $X_1,\ldots,X_n$, let $\hat\bbeta$ be the parameter estimate and $\hat Z_1,\ldots,\hat{Z}_n$ be the estimated residuals.  A set of bootstrapped residuals can be obtained as follows:
\begin{enumerate}
\item
	Let $\hat{F}_n$ be the mean-corrected empirical distribution of $\{\hat Z_j,1\le j \le n\}$;
\item
	Generate $X^*_1,\ldots,X^*_n$ from the time series model with parameter value $\hat\bbeta$ and innovation sequence $\{Z_j^*\}$ generated from $\hat{F}_n$;
\item
	Re-fit the time series model.  Obtained the parameter estimate $\hat\bbeta^*$ and the estimated residuals $\hat{Z}^*_1,\ldots,\hat{Z}^*_n$.
\end{enumerate}
Let $n \hat{T}_h(\hat Z^*,\mu)$ be the ADCV calculated from the bootstrapped residuals $\hat{Z}^*_1,\ldots,\hat{Z}^*_n$.  
In Theorem~\ref{thm:boot} below, we show that when the sample size $n$ is large, the empirical distribution of $\{n\hat{T}_h(\hat Z^*,\mu)\}$ forms a good representation of the limiting distribution of $n \hat{T}_h(\hat Z,\mu)$, the ADCV of the actual fitted residuals. Before stating the theorem, we first state the relevant conditions.  We denote by $\P_n$ and $\E_n$ the probability and expectation conditional on the observations $X_1,\ldots,X_n$.
\begin{enumerate}
\item[(M1')]\label{cond:m1p}	
	Let $\calF_j,\calF_j^*$ be the $\sigma$-algebra generated by $\{Z_k,k\le j\}$ and $\{Z^*_k,k\le j\}$, respectively.  We assume that condition \hyperref[cond:m3]{(M1)} holds, i.e., \eqref{eq:m11} and \eqref{eq:m12} hold.  In addition, as $n\to\infty$, for any $\epsilon>0$,
	$$
		\P_n\left(\left|\frac{1}{n}\sum_{j=1}^n\E_n[\boldm^T(Z_{-\infty:j}^{*};\hat\bbeta)\boldm(Z_{-\infty:j}^{*};\hat\bbeta)|\calF_{j-1}^{*}] - \tau^2\right|>\epsilon \right) \cip 0
	$$
	for some $\tau>0$, and
	$$
		\P_n\left(\frac{1}{n}\sum_{j=1}^n\E_n[\boldm^T(Z_{-\infty:j}^{*};\hat\bbeta)\boldm(Z_{-\infty:j}^{*};\hat\bbeta) {\bf1}_{\{|\boldm(Z_{-\infty:j}^{*};\hat\bbeta)|>\sqrt{n}\epsilon\}}|\calF_{j-1}^{*}]>\epsilon\right) \cip 0.
	$$
	
	\item[(M2')]\label{cond:m2p}
	Assume that the function $h$ in the invertible representation \eqref{eq:invert} is continuously differentiable and
	\beqo\label{eq:3}
		\P\left[\sup_n\E_n\|\bl_0^*(\hat\bbeta)\|^2 <\infty\right] =1,
	\eeqo
	where
	\beqq \label{eq:bigL:b}
		\bl^*_j(\bbeta) := \frac{\partial}{\partial\bbeta} h(X^*_{-\infty:j};\bbeta).
	\eeqq
\item[(M3')]\label{cond:m3p}
	Assume that the estimated residuals based on the finite sequence of observations, $\hat Z^*_{j}$, is close to the fitted residuals based on the infinite sequence, $\tilde Z^*_{j}$, such that for any $\epsilon>0$,
	$$
		\P_n\left(\frac{1}{\sqrt{n}} \sum_{j=1}^n |\hat{Z}^*_{j} - \tilde{Z}^*_{j}|^k >\epsilon\right) \to 0,  \quad k=1,2.
	$$
\end{enumerate}

\bre
Condition~\hyperref[cond:m1p]{(M1')} ensures that $\sqrt{n}(\hat\bbeta^* - \hat\bbeta)$ provides a good approximation to $\bQ$, the limit of $\sqrt{n}(\hat\bbeta - \bbeta_0)$.  These conditions are standard for the martingale central limit theorem, see, for example, \cite{scott:1973}.  Conditions~\hyperref[cond:m2p]{(M2')} and \hyperref[cond:m3p]{(M3')} are parallel arguments to conditions~\hyperref[cond:m2]{(M2)} and \hyperref[cond:m3]{(M3)}.
\ere

\bth \label{thm:boot}
Assuming conditions \hyperref[cond:m1p]{(M1')}, \hyperref[cond:m2p]{(M2')} and \hyperref[cond:m3p]{(M3')} hold, the ADCV of the bootstrapped residuals $\{\hat{Z}^*_{1:n}\}$ satisfies
$$
	\sup_t\left|\P_n\left( n \hat{T}_h(\hat Z^*,\mu)\le t\right) - \P\left(\|G_h + \xi_h\|^2_\mu\le t\right)\right| \cip 0.
$$
\ethe



\section{Example: ARMA($p$,$q$)} \label{sec:arma}

\noindent
Consider the causal, invertible ARMA($p,q$) process that follows the recursion,
\beqq \label{eq:arma}
	X_t = \sum_{i=1}^p \phi_i X_{t-i} + Z_t + \sum_{j=1}^q \theta_j Z_{t-j},
\eeqq
where $\bbeta=(\phi_1,\ldots,\phi_p,\theta_1,\ldots,\theta_q)^T$ is the vector of parameters and $\{Z_t\}$ is iid with mean 0 and variance $\sigma^2$.
Denote the AR and MA polynomials by $\phi(z) = 1-\sum_{i=1}^p \phi_i z^i$ and $\theta(z) = 1+\sum_{j=1}^q \theta_j z^j$, and let $B$ be the backward operator such that
$$
	BX_t = X_{t-1}.
$$
Then the recursion \eqref{eq:arma} can be represented by
$$
	\phi(B) X_t = \theta(B) Z_t.
$$
It follows from invertibility that $\phi(z)/\theta(z)$ has the power series expansion
$$
	\frac{\phi(z)}{\theta(z)} = \sum_{j=0}^\infty \pi_j(\bbeta) z^i,
$$
where $\sum_{j=0}^\infty |\pi_j(\bbeta)| < \infty$, and
$$
	Z_t = Z_t(\bbeta) = \sum_{j=0}^\infty \pi_j(\bbeta) X_{t-j}.
$$
Given an estimate of the parameters $\hat\bbeta$, the residuals based on the infinite sequence $\{X_{-\infty:n}\}$ are given by
$$
	\tilde{Z}_t := Z_t(\hat\bbeta) = \sum_{j=0}^\infty \pi_j(\hat\bbeta)  X_{t-j}.
$$
Based on the observed data $X_1,\ldots, X_n$, the estimated residuals are
\beqq \label{eq:arma:residuals}
	\hat{Z}_t = \sum_{j=0}^{t-1} \pi_j(\hat\bbeta) X_{t-j}.
\eeqq
One choice for $\hat\bbeta$ is the pseudo-MLE based on Gaussian likelihood 
$$
	L(\bbeta,\sigma^2) \propto \sigma^{-n} |\bSig|^{-1/2} \exp\{\frac{1}{2\sigma^2}\bfX_n^T\bSig^{-1}\bfX_n\},
$$
where $\bfX_n = (X_1,\ldots,X_n)^T$ and the covariance $\bSig=\bSig(\bbeta):=\var(\bfX_n)/\sigma^2$ is independent of $\sigma^2$.  The pseudo-MLE $\hat\bbeta$ and $\hat\sigma^2$ are taken to be the values that maximize $L(\bbeta,\sigma^2)$. It can be shown that $\hat\bbeta$ is consistent and asymptotically normal even for non-Gaussian $Z_t$ \citep{brockwell:davis:1991}.  

We have the following result for the ADCV of ARMA residuals.

\bco \label{thm:arma}
Let $\{X_t,1\le j \le n\}$ be observations from a causal and invertible ARMA($p$,$q$) time series and $\{\hat Z_t,1\le t \le n\}$ be the estimated residuals defined in \eqref{eq:arma:residuals} using the pseudo-MLE $\hat\bbeta$. Assume that $\mu$ satisfies \eqref{eq:7}, then
$$
	n \hat{T}_h(\hat Z;\mu) \cid \|G_h + \xi_h\|_\mu^2,
$$
where $(G_h,\xi_h)$ is a joint Gaussian process defined on $\bbr^2$ with $G_h$ as specified in \eqref{eq:gh} and $\xi_h$ in \eqref{eq:xi}.
\eco
The proof of Corollary~\ref{thm:arma} is given in Appendix~\ref{app:arma}.

\bre
In the case where the distribution of $Z_t$ is in the domain of attraction of an $\alpha$-stable law with $\alpha\in (0,2)$, and the parameter estimator $\hat\bbeta$ has convergence rate faster than $n^{-1/2}$, i.e.,
$$
	a_n(\hat\bbeta-\bbeta) = O_p(1), \quad \text{for some $a_n=o(n^{-1/2})$},
$$
\citep{davis:1996}, the ADCV of the residuals has limit
$$
	n \hat{T}_h(\hat Z;\mu) \cid \|G_h \|_\mu^2,
$$
where the correction term $\xi_h$ disappears.  For a proof in the AR($p$) case, see Theorem~4.2 of \cite{davis:matsui:mikosch:wan:2018}.
\ere

\subsection{Simulation}
\begin{figure}[ht]
	\centering
	\includegraphics[width=5in]{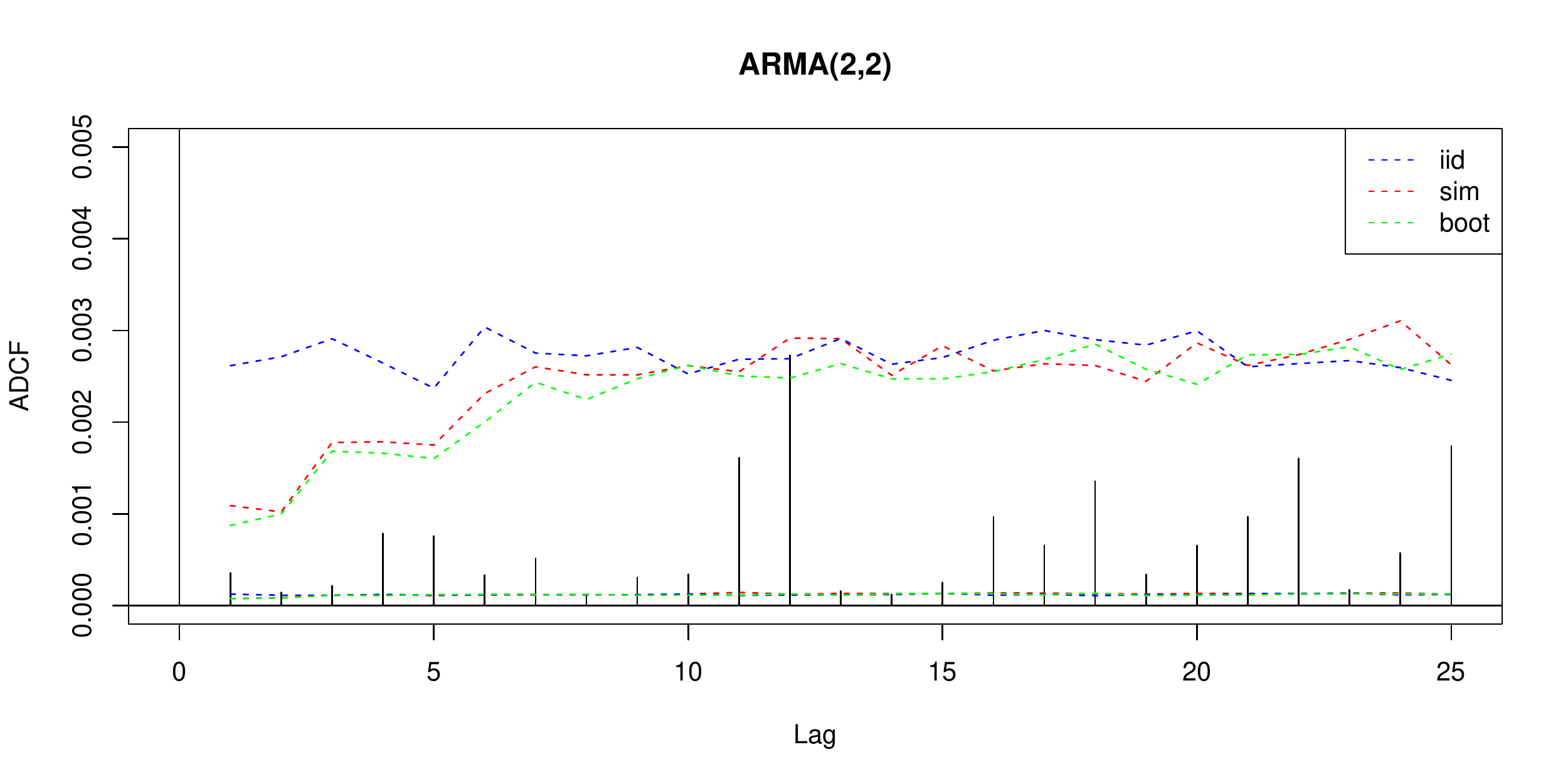}
	\caption{Empirical $5\%$ and $95\%$ quantiles of the ADCF for a) iid innovations; b) estimated residuals; c) bootstrapped residuals; from a ARMA(2,2) model. }
	\label{fig:arma}
\end{figure}

We generate time series of length $n=2000$ from an ARMA(2,2) model with standard normal innovations and parameter values 
$$
	\bbeta=(\phi_1,\phi_2,\theta_1,\theta_2)=(1.2, -0.32, -0.2, -0.48).
$$
For each simulation, an ARMA(2,2) model is fitted to the data.  In Figure~\ref{fig:arma}, we compare the empirical $5\%$ and $95\%$ quantiles for the ADCF of
\begin{enumerate}
\item[a)]
	iid innovations from 1000 independent simulations;
\item[b)]
	estimated residuals from 1000 independent simulations of $\{X_t\}$;
\item[c)]
	estimated residuals through 1000 independent parametric bootstrap samples from one realization of $\{X_t\}$.
\end{enumerate}
In order to satisfy condition \eqref{eq:7}, the ADCFs are evaluated using the Gaussian weight measure $N(0,0.5^2)$.
Confirming the results in Theorem~\ref{thm:meta} and Corollary~\ref{thm:arma}, the simulated quantiles of $\hat{R}_h(\hat Z;\mu)$ differ significantly from that of $\hat{R}_h(Z;\mu)$, especially when $h$ is small.  Given one realization of the time series, the quantiles estimated by parametric boostrap correctly capture this effect.

\section{Example: GARCH($p$,$q$)} \label{sec:garch}

In this section, we consider the GARCH($p$,$q$) model,
$$
	X_t = \sigma_t Z_t,
$$
where the $Z_t$'s are iid innovations with mean 0 and variance 1 and
\beqq \label{eq:cv:def}
	\sigma_t^2 = \alpha_0 + \sum_{i=1}^p \alpha_i X_{t-i}^2 + \sum_{j=1}^q \beta_j \sigma_{t-j}^2,\quad \alpha_0>0,\, \alpha_i \ge 0,\, \beta_j\ge0.
\eeqq
Let $\btheta=(\alpha_0,\alpha_1,\ldots,\alpha_p,\beta_1,\ldots,\beta_q)$ denote the parameter vector.  We write the conditional variance $\sigma_t^2=\sigma_t^2(\btheta)$ to denote it as a function of $\btheta$.

Iterating the recursion in \eqref{eq:cv:def} gives 
$$
	\sigma_t^2(\btheta) = c_0(\btheta) + \sum_{i=1}^\infty c_i(\btheta) X_{t-i}^2,
$$
for suitably defined functions $c_i$'s, see \cite{berkes:horvath:kokoszka:2003}.  
Given an estimator $\hat\btheta$, an estimator for $\sigma^2_t(\btheta)$ based on the infinite sequence $\{X_j,j\le t\}$ can be written as
$$
	\tilde\sigma_t^2 := \sigma_t^2(\hat\btheta_n) =c_0(\hat\btheta_n) + \sum_{i=1}^\infty c_i(\hat\btheta_n) X_{t-i}^2,
$$
and the unobserved residuals are given by
$$
	\tilde{Z}_t = X_t / \tilde\sigma_t.
$$
In practice, $\tilde\sigma_t^2$ can be approximated by the truncated version
$$	
	\hat\sigma_t^2(\hat\btheta_n) :=c_0(\hat\btheta_n) + \sum_{i=1}^{t-1} c_i(\hat\btheta_n) X_{t-i}^2,
$$
and the estimated residual $\hat Z_t$ is given by
\beqq \label{eq:garch:residuals}
\hat{Z}_t = X_t / \hat\sigma_t.
\eeqq

Define the parameter space by
$$
	\bTheta = \{\bu=(s_0,s_1,\ldots,s_p,t_1,\ldots,t_q): t_1+\cdots+t_q\le \rho_0, \underline{u} \le \min(\bu) \le \max(\bu) \le \bar{u}\},
$$
for some $0<\underline{u}<\bar{u}$, $0<\rho_0<1$ and $q\underline{u}<\rho_0$, and assume the following conditions:
\begin{enumerate}
\item[(Q1)] \label{cond:q1}
	The true value $\btheta$ lies in the interior of $\bTheta$.
\item[(Q2)] \label{cond:q2}
	For some $\zeta>0$,
	$$\lim_{x\to0} x^{-\zeta}\P\{|Z_0|\le x\} = 0.$$
\item[(Q3)] \label{cond:q3}
	For some $\delta>0$,
	$$\E|Z_0|^{4+\delta} < \infty.$$
\item[(Q4)] \label{cond:q4}
	The GARCH($p,q$) representation is minimal, i.e., the polynomials $A(z)= \sum_{i=1}^p\alpha_iz^i$ and $B(z)=1-\sum_{j=1}^p\beta_jz^j$ do not have common roots.
\end{enumerate}
Given observations $\{X_t, 1\le t\le n\}$, \cite{berkes:horvath:kokoszka:2003} proposed a quasi-maximum likelihood estimator for $\btheta$ given by
$$
	\hat\btheta_n := {\arg\max}_{\bu\in\bTheta} \sum_{t=1}^n l_t(\bu),
$$
where
$$
	l_t(\bu) := -\frac12\log \hat\sigma_t^2(\bu) - \frac{X_t^2}{2\hat\sigma_t^2(\bu)}.
$$
Provided that \hyperref[cond:q1]{(Q1)}--\hyperref[cond:q4]{(Q4)} are satisfied, the quasi-MLE $\hat\btheta_n$ is consistent and asymptotically normal.

Consider the estimated residuals for the GARCH($p$,$q$) model based on $\hat\btheta_n$. We have the following result.

\bco \label{thm:garch}

Let $\{X_t,1\le j \le n\}$ be observations from a GARCH($p$,$q$) time series and $\{\hat Z_t,1\le t \le n\}$ be the estimated residuals defined in \eqref{eq:garch:residuals} based on the quasi-MLE $\hat\btheta_n$.  Assume that \hyperref[cond:q1]{(Q1)}--\hyperref[cond:q4]{(Q4)} holds and that $\mu$ satisfies \eqref{eq:7},
we have
$$
	n \hat{T}_h(\hat Z;\mu) \cid \|G_h + \xi_h\|_\mu^2,
$$
where $(G_h,\xi_h)$ is a joint Gaussian process defined on $\bbr^2$ with $G_h$ as specified in \eqref{eq:gh} and $\xi_h$ in \eqref{eq:xi}.
\eco
The proof of Corollary~\ref{thm:garch} is given in Appendix~\ref{app:garch}.

\subsection{Simulation}

\begin{figure}[ht]
	\centering
	\includegraphics[width=5in]{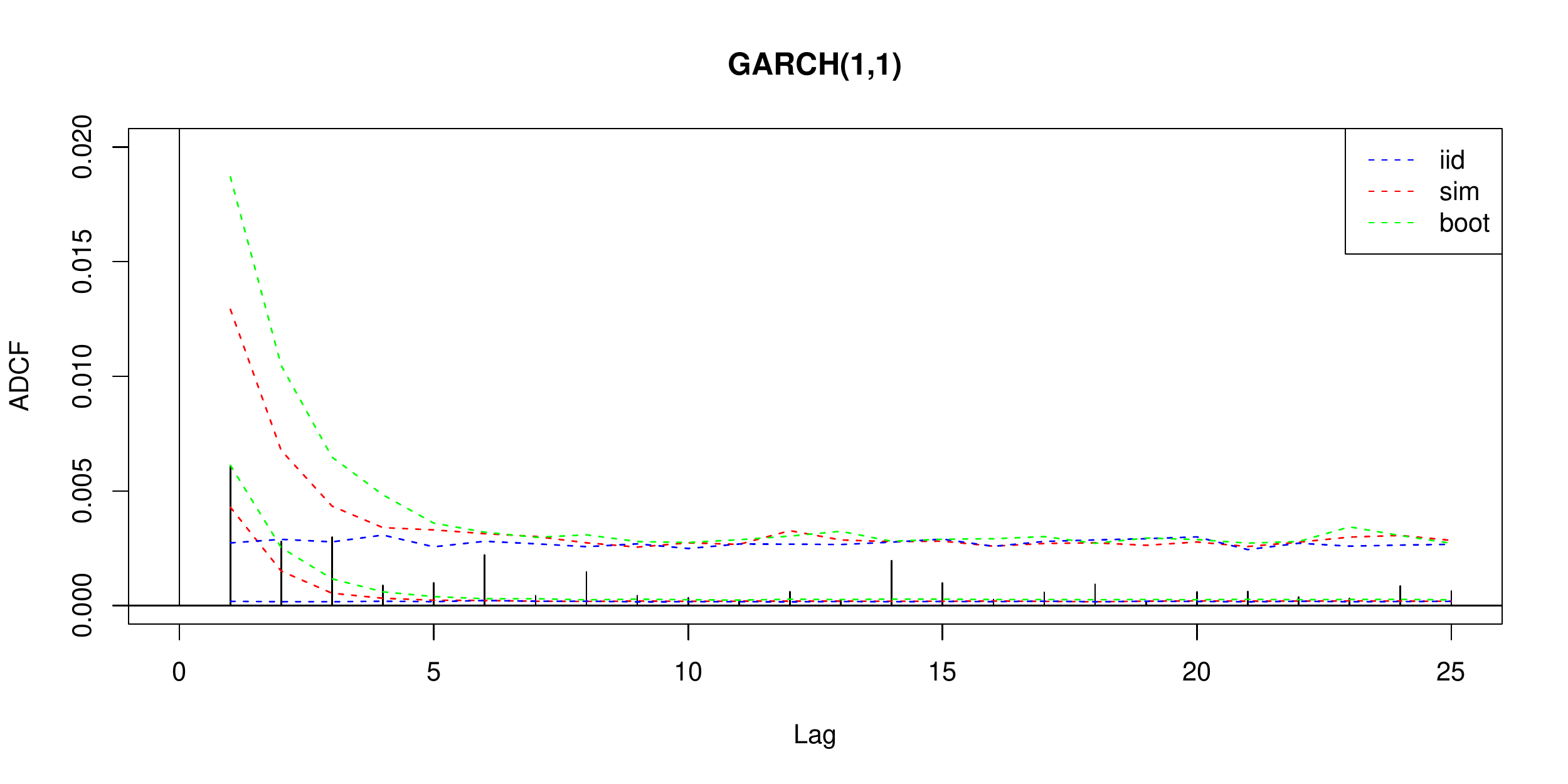}
	\caption{Empirical $5\%$ and $95\%$ quantiles of the ADCF for a) iid innovations; b) estimated residuals; c) bootstrapped residuals; from a GARCH(1,1) model. }
	\label{fig:garch}
\end{figure}

We generate time series of length $n=2000$ from a GARCH(1,1) model with parameter values 
$$
	\btheta=(\alpha_0,\alpha_1,\beta_1)=(0.5,0.1,0.8).
$$
For each simulation, a GARCH(1,1) model is fitted to the data.
In Figure~\ref{fig:garch}, we compare the empirical $5\%$ and $95\%$ quantiles for the ADCF of
\begin{enumerate}
\item[a)]
	iid innovations from 1000 independent simulations;
\item[b)]
	estimated residuals from 1000 independent simulations of $\{X_t\}$;
\item[c)]
	estimated residuals through 1000 independent parametric bootstrap samples from one realization of $\{X_t\}$.
\end{enumerate}
Again the ADCFs are based on the Gaussian weight measure $N(0,0.5^2)$. The difference between the quantiles of $\hat{R}_h(\hat Z;\mu)$ and $\hat{R}_h(Z;\mu)$ can be observed.  For this GARCH model, the correction has the opposite effect than in the previous ARMA exaple -- the ADCF for residuals are larger than that for iid variables, especially for small lags.

\section{Example: Non-causal AR(1)} \label{sec:noncausal}

\begin{figure}[ht]
	\centering
	\includegraphics[width=5in]{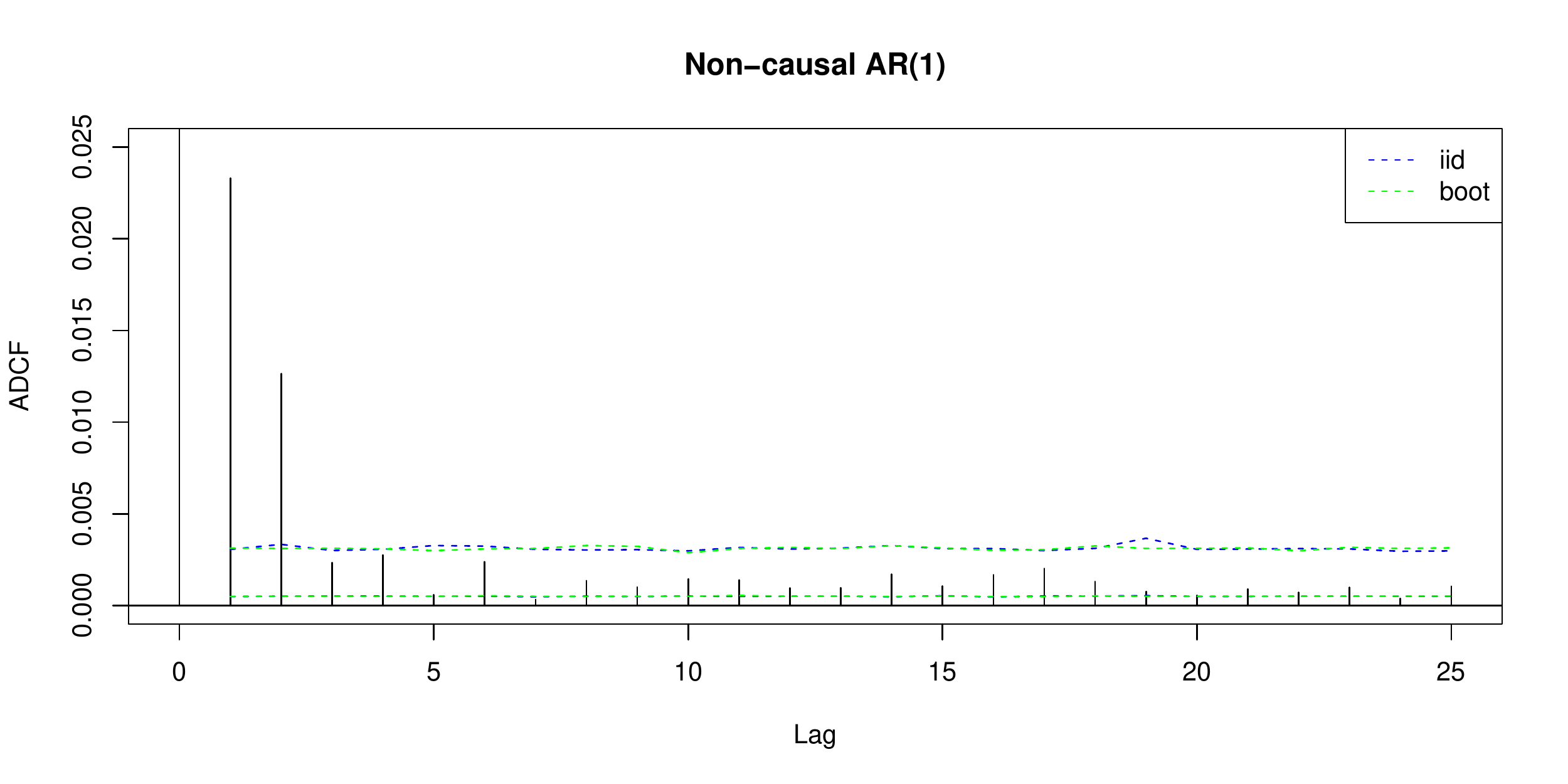}
	\caption{Empirical $5\%$ and $95\%$ quantiles of the ADCF for a) iid innovations; b) bootstrapped residuals; from non-causal AR(1) data fitted with a causal AR(1) model. }
	\label{fig:ncar}
\end{figure}

In this section, we consider an example where the model is misspecified.  We generate time series of length $n=2000$ from a non-causal AR(1) model 
$$
	X_t = \phi X_{t-1}+Z_t
$$
with $\phi=1.67$ and $Z_t$'s from a $t$-distribution with 2.5 degrees of freedom.  Then we fit a causal AR(1) model, where $|\phi|<1$, to the data and obtain the corresponding residuals. Again we use the Gaussian weight measure $N(0,0.5^2)$ when evaluating the ADCF of the residuals. In Figure~\ref{fig:ncar}, the $5\%$ and $95\%$ ADCF quantiles are plotted for:
\begin{enumerate}
\item[a)]	
	estimated residuals from 1000 independent simulations of $\{X_t\}$;
\item[b)]
	estimated residuals through 1000 independent parametric bootstrap samples from one realization of $\{X_t\}$.
\end{enumerate}

The ADCFs of the bootstrapped residuals provide an approximation for the limiting distribution of the ADCF of the residuals given the model is correctly specified.  In this case, the ADCFs of the estimated residuals significantly differ from the quantiles of that of the bootstrapped residuals.  This indicates the time series does not come from the assumed causal AR model.

\section{Conclusion}

In this paper, we propose a goodness-of-fit procedure for time series models by examining the serial dependence of estimated residuals.  The dependence is measured using the auto-distance covariance function (ADCV) and its limiting behavior is derived for general classes of time series models.  We show that the limiting law often differs from that of the ADCV based on iid innovations by a correction term.  This indicates that adjustments should be made when testing the goodness-of-fit of the model.  We illustrate the result on simulated examples of ARMA and GARCH processes and discover that the adjustments could be in either direction -- the quantiles of ADCV for residuals could be larger or smaller than that for iid innovations.  We also studied an example when a non-causal AR process was incorrectly fitted with a causal model and showed that ADCV correctly detected model misspecification when applied to the residuals.

\bibliographystyle{plainnat}

\bibliography{wan}

\newpage
\appendix

In the following appendices, we provide proofs to Theorem~\ref{thm:meta} and Corollaries~\ref{thm:arma} and \ref{thm:garch}.
Throughout the proofs, $c$ denotes a general constant whose value may change from line to line.

\section{Proof of Theorem~\ref{thm:meta}}\label{app:meta}
\bpf

The proof proceeds in the following steps with the aid of Propositions~\ref{prop:joint:conv}, \ref{prop:a2} and \ref{prop:a3}.
Write
$$
n \hat{T}_h(\hat Z;\mu) =: \|\sqrt{n}C_n^{\hat Z}\|^2_\mu = \|\sqrt{n}C_n^{\hat Z} - \sqrt{n}C_n^{Z} + \sqrt{n}C_n^{Z}\|^2_\mu,
$$
where
$$
	C_n^{\hat{Z}}(s,t) := \frac{1}{n} \sum_{j=1}^{n-h} e^{is\hat{Z}_{j}+it\hat{Z}_{j+h}} - \frac{1}{n} \sum_{j=1}^{n-h} e^{is\hat{Z}_{j}} \frac{1}{n} \sum_{j=1}^{n-h} e^{it\hat{Z}_{j+h}}
$$
and
$$
	C_n^{Z}(s,t) := \frac{1}{n} \sum_{j=1}^{n-h} e^{isZ_{j}+itZ_{j+h}} - \frac{1}{n} \sum_{j=1}^{n-h} e^{isZ_{j}} \frac{1}{n} \sum_{j=1}^{n-h} e^{itZ_{j+h}}.
$$
We first show in Proposition~\ref{prop:joint:conv} that
\beqo \label{eq:jointconv}
	(\sqrt{n}(C_n^{\hat Z} - C_n^{Z}), \sqrt{n}C_n^{Z}) \cid (\xi_h,G_h), \quad \text{on $\mathcal{C}(K)$},
\eeqo
where $K$ is any compact set in $\bbr^2$.  This implies
\beqo \label{eq:convsum}
	\sqrt{n}C_n^{\hat Z} \cid \xi_h+G_h, \quad \text{on $\mathcal{C}(K)$}.
\eeqo
For $\delta\in(0,1)$, define the compact set
$$
	K_\delta = \{(s,t)|\delta\le s \le 1/\delta,\,\delta\le t\le 1/\delta \}.
$$
It follows from the continuous mapping theorem that
$$
	n\int_{K_\delta} |C_n^{\hat Z} |^2 \mu(ds,dt) \cid \int_{K_\delta} |G_h + \xi_h |^2 \mu(ds,dt).
$$
To complete the proof, it remains to justify that we can take $\delta\downarrow0$.  For this it suffices to show that for any $\varepsilon>0$,
$$
	\lim_{\delta\to0} \limsup_{n\to\infty} \P\left( \int_{K_\delta^c}|\sqrt{n}C_n^{\hat Z}|^2 \mu(ds,dt)>\varepsilon \right) = 0,
$$
and
$$
	\lim_{\delta\to0} \P\left( \int_{K_\delta^c} |G_h + \xi_h |^2 \mu(ds,dt)>\varepsilon \right) = 0.
$$
These are shown in Propositions~\ref{prop:a2} and \ref{prop:a3}, respectively.

\epf

\bpr \label{prop:joint:conv}
Given the conditions \hyperref[cond:m1]{(M1)}--\hyperref[cond:m3]{(M3)},
$$
	(\sqrt{n}(C_n^{\hat Z} - C_n^{Z}), \sqrt{n}C_n^{Z}) \cid (\xi_h,G_h), \quad \text{on $\mathcal{C}(K)$},
$$
for any compact $K\subset \bbr^2$.

\epr

\bpf
We first consider the marginal convergence of $\sqrt{n}(C_n^{\hat Z} - C_n^{Z})$.  Denote
$$
	E_n(s,t) :=  \frac{1}{\sqrt{n}} \sum_{j=1}^{n-h} \left( \ex^{is\hat Z_j+it\hat Z_{j+h}} -  \ex^{isZ_j+itZ_{j+h}} \right),
$$
then
\beam
	\sqrt{n}(C_n^{\hat Z}(s,t) - C_n^{Z}(s,t) )
		&=& \frac{1}{\sqrt{n}} \sum_{j=1}^{n-h} \left( \ex^{is\hat Z_j+it\hat Z_{j+h}} -  \ex^{isZ_j+itZ_{j+h}} \right)\nonumber\\
		&& \quad - \frac{1}{\sqrt{n}} \sum_{j=1}^{n-h} \left(\ex^{is\hat Z_j}  - \ex^{isZ_j} \right)\frac{1}{n}\sum_{j=1}^{n-h} \ex^{itZ_{j+h}}\nonumber\\
		&& \qquad - \frac{1}{n} \sum_{j=1}^{n-h} \ex^{is\hat Z_j} \frac{1}{\sqrt{n}}\sum_{j=1}^{n-h} \left(\ex^{it\hat Z_{j+h}} - \ex^{itZ_{j+h}} \right)\nonumber\\
		&=& E_n(s,t) - E_n(s,0)\frac{1}{n}\sum_{j=1}^{n-h} \ex^{itZ_{j+h}} - E_n(0,t)\frac{1}{n} \sum_{j=1}^{n-h} \ex^{is\hat Z_j}. \label{eq:decomp:en}
\eeam
We now derive the limit of $E_{n}(s,t)$.
Observe that uniformaly for $(s,t)\in K$,
\beao
	E_{n}(s,t)
	&=& \frac{1}{\sqrt{n}} \sum_{j=1}^{n-h} \ex^{isZ_j+itZ_{j+h}} \left(\ex^{is(\hat Z_j - Z_j)+it(\hat Z_{j+h}-Z_{j+h})} -1 \right) \\
		&=& \frac{1}{n} \sum_{j=1}^{n-h} \ex^{isZ_j+itZ_{j+h}} (is\sqrt{n}(\hat Z_j - Z_j)+it\sqrt{n}(\hat Z_{j+h}-Z_{j+h})) + o_p(1),\\
		&=& \frac{1}{n} \sum_{j=1}^{n-h} \ex^{isZ_j+itZ_{j+h}} (is\sqrt{n}(\hat Z_j -  \tilde Z_j)+it\sqrt{n}(\hat Z_{j+h}- \tilde Z_{j+h})) \\
		&& \quad + \frac{1}{n} \sum_{j=1}^{n-h} \ex^{isZ_j+itZ_{j+h}} (is\sqrt{n}(\tilde Z_j - Z_j)+it\sqrt{n}(\tilde Z_{j+h}-Z_{j+h})) + o_p(1) \nonumber\\
		&=:& E_{n1}(s,t) + E_{n2}(s,t) + o_p(1).
\eeao
By assumption \hyperref[cond:m3]{(M3)},
$$
	|E_{n1}(s,t)| \le |s| \frac{1}{\sqrt{n}} \sum_{j=1}^{n-h} |\hat Z_j -  \tilde Z_j| + |t| \frac{1}{\sqrt{n}} \sum_{j=1}^{n-h} |\hat Z_{j+h} -  \tilde Z_{j+h}| \cip 0, \quad \text{in $\mathcal{C}(K)$.}
$$
It follows from a Taylor expansion that
\beao
		E_{n2}(s,t) &=& \sqrt{n}(\hat\bbeta - \bbeta)^T \frac{1}{n} \sum_{j=1}^{n-h} \ex^{isZ_j+itZ_{j+h}} \left(is \bl_j(\bbeta^*)+it \bl_{j+h}(\bbeta^*)\right),
\eeao
where $\bbeta^* = \bbeta + \epsilon(\hat\bbeta-\bbeta)$ for some $\epsilon\in[0,1]$.
Since $\bl_j(\bbeta)$ is stationary and ergodic, 
in view of the uniform ergodic theorem,
$$
	\frac{1}{n} \sum_{j=1}^{n-h} \ex^{isZ_j+itZ_{j+h}} \left(is\bl_j(\bbeta)+it \bl_{j+h}(\bbeta)\right) 
	\cip \E\left[\ex^{isZ_j+itZ_{j+h}} \left(is \bl_j(\bbeta)+it \bl_{j+h}(\bbeta)\right)\right] =: \mathbf C_h(s,t), \quad \text{in $\mathcal{C}(K)$.}
$$
Hence,
$$
	E_{n}(s,t) \cid \bQ^T\mathbf C_h(s,t), \quad \text{in $\mathcal{C}(K)$.}
$$
Note that
$$
	\frac{1}{n}\sum_{j=1}^{n-h} \ex^{itZ_{j+h}} \cip \varphi_Z(t), \quad \text{in $\mathcal{C}(K)$,}
$$
and
\beqq\label{eq:conv:res}
	\frac{1}{n} \sum_{j=1}^{n-h} \ex^{is\hat Z_j} = \frac{1}{n} \sum_{j=1}^{n-h} \ex^{is Z_j} + \frac{1}{\sqrt{n}} E_{n}(s,0)  \cip \varphi_Z(s), \quad \text{in $\mathcal{C}(K)$.}
\eeqq
We have
$$
	\sqrt{n}(C_n^{\hat Z} - C_n^Z) \cid \bQ^T \left(\mathbf C_h(s,t) - \mathbf C_h(s,0)\varphi_Z(t) - \mathbf C_h(0,t) \varphi_Z(s) \right), \quad \text{in $\mathcal{C}(K)$.}
$$
To further simplify the above expression, notice that $\bl_j(\bbeta)$ is a function of $X_{-\infty:j}$ and independent of $Z_{j+h}$ by causality.  Hence
\beao
	\mathbf C_h(s,t) &=& \E\left[\ex^{isZ_j} is \bl_j(\bbeta)\right] \E\left[ \ex^{itZ_{j+h}}\right]+ \E\left[\ex^{isZ_j+itZ_{j+h}} it \bl_{j+h}(\bbeta)\right] \\
	&=& \mathbf C_h(s,0) \varphi_Z(t) + \E\left[\ex^{isZ_j+itZ_{j+h}} it \bl_{j+h}(\bbeta)\right],
\eeao
and
\beam
	 &&\hspace{-.5in}\bQ^T \left(\mathbf C_h(s,t) - \mathbf C_h(s,0)\varphi_Z(t) - \mathbf C_h(0,t) \varphi_Z(s)\right) \nonumber\\
	  &=& \bQ^T \left(\E\left[\ex^{isZ_j+itZ_{j+h}} it \bl_{j+h}(\bbeta)\right] - \E\left[\ex^{itZ_{j+h}} it \bl_{j+h}(\bbeta)\right] \varphi_Z(s)\right) 
	\,=\, \xi_h(s,t). \label{eq:xi:cal}
\eeam
This justifies the marginal convergence of $\sqrt{n}(C_n^{\hat Z} - C_n^Z)$. 

For the joint convergence of $\sqrt{n}(C_n^{\hat Z} - C_n^Z)$ and $\sqrt{n}C_n^Z$, we recall assumption \hyperref[cond:m1]{(M1)}
$$
	\sqrt{n}(\hat\bbeta-\bbeta) = \frac{1}{\sqrt{n}} \sum_{j=1}^{n} \boldm(X_{-\infty:j};\bbeta) + o_p(1) 
$$
and also note from the proof of Theorem~1 in \cite{davis:matsui:mikosch:wan:2018} that
$$
	\sqrt{n}C_n^Z = \frac{1}{\sqrt{n}} \sum_{j=1}^{n} (e^{isZ_j}-\varphi_Z(s)) (e^{itZ_{j+h}}-\varphi_Z(t)) + o_p(1)\cid G_h, \quad\text{in $\mathcal{C}(K)$.}
$$
By martingale central limit theorem,
$$
	\left(\frac{1}{\sqrt{n}} \sum_{j=1}^{n} \boldm(X_{-\infty:j};\bbeta),\quad \frac{1}{\sqrt{n}} \sum_{j=1}^{n-h} (e^{isZ_j}-\varphi_Z(s)) (e^{itZ_{j+h}}-\varphi_Z(t))\right)
$$
converges jointly to $(\bQ,G_h)$.
This implies the joint convergence of $\sqrt{n}(\hat\bbeta-\bbeta)$ and $\sqrt{n}C_n^Z$.
Since $\xi_h$ continuous and its randomness only depends on $\bQ$, the joint convergence $\sqrt{n}C_n^Z$ and $\sqrt{n}C_n^{\hat{Z}} - \sqrt{n}C_n^Z$ also follows.

\epf

\bpr \label{prop:a2}
Under the conditions of Theorem~3.1,
$$
	\lim_{\delta\to0} \limsup_{n\to\infty} \P\left( \int_{K_\delta^c}|\sqrt{n}C_n^{\hat Z}|^2 \mu(ds,dt)>\varepsilon \right) = 0.
$$
\epr

\bpf
Using telescoping sums, $C_n^{\hat Z} - C_n^{Z}$ has the following decomposition,
\beao
C_n^{\hat Z} - C_n^{Z} 
	&=& \frac{1}{n} \sum_{j=1}^{n-h} A_jB_j - \frac{1}{n} \sum_{j=1}^{n-h} A_j\frac{1}{n} \sum_{j=1}^{n-h} B_j - \frac{1}{n} \sum_{j=1}^{n-h} U_j\frac{1}{n} \sum_{j=1}^{n-h} B_j  - \frac{1}{n} \sum_{j=1}^{n-h} V_j\frac{1}{n} \sum_{j=1}^{n-h} A_j \\
	&&\quad +  \frac{1}{n} \sum_{j=1}^{n-h} U_jB_j  + \frac{1}{n}  \sum_{j=1}^{n-h} V_jA_j \,=:\, \sum_{k=1}^6 I_{nk}(s,t),
\eeao
where
\beao
	U_j = e^{isZ_j} - \varphi_Z(s),\quad V_j =  e^{itZ_{j+h}} - \varphi_Z(t), \quad A_j = e^{is\hat{Z}_j} - e^{isZ_j}, \quad B_j = e^{it\hat{Z}_{j+h}} - e^{itZ_{j+h}}.
\eeao
From a Taylor expansion,
\beao
	n|I_{n1}(s,t)|^2 
	&\le& \left( \frac{1}{\sqrt{n}} \sum_{j=1}^{n-h} |A_jB_j|\right)^2 \\
	&\le& \left( \frac{1}{\sqrt{n}} \sum_{j=1}^{n-h} |e^{is(\hat{Z}_j-Z_j)}-1||e^{it(\hat{Z}_{j+h}-Z_{j+h})}-1|\right)^2 \\
	&\le& c\ \left( \frac{1}{\sqrt{n}} \sum_{j=1}^{n-h} \left(1\wedge|s||\hat{Z}_j-Z_j|\right) \left(1\wedge|t||\hat{Z}_{j+h}-Z_{j+h}|\right) \right)^2 \\
	&\le& c\ \min \left( |s|^2  \left(\frac{1}{\sqrt{n}}\sum_{j=1}^{n-h} |\hat{Z}_j-Z_j|\right)^2,\, |t|^2 \left(\frac{1}{\sqrt{n}}\sum_{j=1}^{n-h} |\hat{Z}_{j+h}-Z_{j+h}|\right)^2,\, \right. \\
	&& \quad \left. |st|^2 \left(\frac{1}{\sqrt{n}}\sum_{j=1}^{n-h} |\hat{Z}_j-Z_j||\hat{Z}_{j+h}-Z_{j+h}|\right)^2\right) \\
	&\le& c\ \min \left( |s|^2  \left(\frac{1}{\sqrt{n}}\sum_{j=1}^{n-h} |\hat{Z}_j-Z_j|\right)^2,\, |t|^2 \left(\frac{1}{\sqrt{n}}\sum_{j=1}^{n-h} |\hat{Z}_{j+h}-Z_{j+h}|\right)^2,\, \right. \\
	&& \quad \left. |st|^2 \left(\frac{1}{\sqrt{n}}\sum_{j=1}^{n-h} |\hat{Z}_j-Z_j|^2\frac{1}{\sqrt{n}}\sum_{j=1}^{n-h}|\hat{Z}_{j+h}-Z_{j+h}|^2\right)\right) 
\eeao
For $k=1,2$,
\beao
	\frac{1}{\sqrt{n}}\sum_{j=1}^{n-h} |\hat{Z}_j-Z_j|^k 
	&\le& c\ \left(\frac{1}{\sqrt{n}}\sum_{j=1}^{n-h} |\hat{Z}_j-\tilde Z_j|^k + \frac{1}{\sqrt{n}}\sum_{j=1}^{n-h} |\tilde{Z}_j-Z_j|^k \right) \\
	&\le& o_p(1) + c\ \frac1{n^{(k-1)/2}} \|\sqrt{n}(\hat\bbeta-\bbeta)\|^k\frac{1}{n}\sum_{j=1}^{n-h}  \|\bl_j(\bbeta^*)\|^k \\
	&=& O_p(1).
\eeao
Therefore
$$
	n|I_{n1}(s,t)|^2 	\le \min(|s|^2,|t|^2,|st|^2) O_p(1) \le \big((1\wedge |s|^2)\,(1\wedge |t|^2) +  (s^2+t^2)\,\1(|s| \wedge |t|>1)\big) O_p(1),
$$
where the $O_p(1)$ term does not depend on $(s,t)$. This implies that
$$
	\lim_{\delta\to0} \limsup_{n\to\infty} \P\left( \int_{K_\delta^c}n|I_{n1}(s,t)|^2 \mu(ds,dt)>\varepsilon \right) = 0.
$$
Similar arguments show that $n|I_{n2}(s,t)|^2$ is bounded by $\min(|s|^2,|t|^2,|st|^2) O_p(1)$, $n|I_{n3}(s,t)|^2$ and $n|I_{n5}(s,t)|^2$ are bounded by $\min(|t|^2,|st|^2) O_p(1)$, and $n|I_{n4}(s,t)|^2$ and $n|I_{n6}(s,t)|^2$ are bounded by $\min(|s|^2,|st|^2) O_p(1)$, and the result of the proposition follows.
\epf

\bpr \label{prop:a3}
Under the conditions of Theorem~3.1,
$$
	\lim_{\delta\to0} \P\left( \int_{K_\delta^c} |G_h + \xi_h |^2 \mu(ds,dt)>\varepsilon \right) = 0.
$$
\epr

\bpf
Note that
\beao
	|\xi(s,t)|^2 
	&\le& c |t|^2 \|\bQ\|^2 \E \left|e^{isZ_0}-\varphi_Z(s)\right|^2 \E|\bl_h(\bbeta)|^2 \\
	&\le& c |t|^2 \|\bQ\|^2 \E \left[\left(1\wedge|s|^2\right)\left(Z_0+\E|Z|\right)^2 \right] \E|\bl_h(\bbeta)|^2 \\
	&\le& |t|^2 \left(1\wedge|s|^2\right) O_p(1).
\eeao
This implies
$$
	\lim_{\delta\to0}  \P\left( \int_{K_\delta^c} |\xi_h |^2 \mu(ds,dt)>\varepsilon \right) = 0.
$$
On the other hand, it was shown in \cite{davis:matsui:mikosch:wan:2018} that $\int|G_h|^2\mu(ds,dt)$ exists as the limit of $n\hat{T}_h(Z;\mu)$.  Hence
$$
	\lim_{\delta\to0}  \P\left( \int_{K_\delta^c} |G_h |^2 \mu(ds,dt)>\varepsilon \right) = 0,
$$
and the proposition is proved.
\epf


\section{Proof of bootstrap consistency: A generalized theorem for triangular arrays}
In this section, we generalize the convergence of ADCV for residuals for triangular arrays, from which the companion result for the bootstrap estimator in Theorem~\ref{thm:boot} can be derived.

Let $\{Z_{1:n}^{(n)}\}$ be a triangular array of random variables where for each $n$, $Z_{j}^{(n)}$'s are defined on the probability space $(\Omega,\calF^{(n)},\P_n)$ such that
$$
	Z_{j}^{(n)} \overset{iid}\sim F^{(n)}.
$$
Assume that the distribution $F^{(n)}$ converges to $F$ in distribution,
$$
	F^{(n)}\cid F,
$$	
where $F$ is the distribution of $Z$.
Let $\{\bbeta^{(n)}\}$ be a sequence of parameter vectors such that
$$
	\bbeta^{(n)} \to \bbeta.
$$
For each $n$, let $\{X_{1:n}^{(n)}\}$ be a time series generated from the model \eqref{eq:model} with parameter vector $\bbeta^{(n)}$ and innovation sequence $\{Z_{-\infty:n}^{(n)}\}$,
$$
	X_{j}^{(n)} = g(Z_{-\infty:j}^{(n)};\bbeta^{(n)}).
$$ 
Let $\hat\bbeta^{(n)}$ and $\{\hat{Z}_{1:n}^{(n)}\}$ be the corresponding estimates and residuals based on $\{X_{1:n}^{(n)}\}$.  We consider $T_n^{(n)}(h)$, the ADCV of $\{\hat{Z}_{1:n}^{(n)}\}$ at lag $h$.  We require the following conditions.

\begin{enumerate}
\item[(N1)]\label{cond:n1}
	Let $\calF_{j}^{(n)}$ be the $\sigma$-algebra generated by $\{Z_{k}^{(n)},k\le j\}$, respectively.  We assume that for any $\epsilon>0$,
	$$
		\P_n\left(\left|\sqrt{n}(\hat\bbeta^{(n)} - \bbeta^{(n)}) - \frac{1}{\sqrt{n}} \sum_{j=1}^n \boldm(Z_{-\infty:j}^{(n)};\bbeta^{(n)}) \right|>\epsilon\right) \to 0.
	$$
	where 
	$$
		\E_n[\boldm(Z_{-\infty:j}^{(n)};\bbeta^{(n)})|\calF_{j-1}^{(n)}] = \mathbf0.
	$$
	Further we assume that as $n\to\infty$, 
	$$
		\frac{1}{n}\sum_{j=1}^n\E_n[\boldm^T(Z_{-\infty:j}^{(n)};\bbeta^{(n)})\boldm(Z_{-\infty:j}^{(n)};\bbeta^{(n)})|\calF_{j-1}^{(n)}] \cip \tau^2,
	$$
	for some $\tau>0$, and
	$$
		\P_n\left(\frac{1}{n}\sum_{j=1}^n\E_n[\boldm^T(Z_{-\infty:j}^{(n)};\bbeta^{(n)})\boldm(Z_{-\infty:j}^{(n)};\bbeta^{(n)}) {\bf1}_{\{|\boldm(Z_{-\infty:j}^{(n)};\bbeta^{(n)})|>\sqrt{n}\epsilon\}}|\calF_{j-1}^{(n)}] >\epsilon\right) \to 0, \quad \forall \epsilon>0.
	$$
\item[(N2)]\label{cond:n2}
	Assume that the function $h$ in the invertible representation \eqref{eq:invert} is continuously differentiable, and writing
	\beqq \label{eq:bigL:b}
		\bl^{(n)}_j(\bbeta) := \frac{\partial}{\partial\bbeta} h(X_{-\infty:j}^{(n)};\bbeta),
	\eeqq
	we have
	\beqo\label{eq:3}
		\sup_n\E_n\|\bl_0^{(n)}(\bbeta^{(n)})\|^2 <\infty.
	\eeqo
\item[(N3)]\label{cond:n3}
	For fixed $j$, let $\tilde Z_{j}^{(n)}$ be the fitted residual based on the unobserved infinite sequence $\{X_{-\infty:j}^{(n)}\}$ obtained from \eqref{eq:tildez}, and $\hat Z_{j}^{(n)}$ be the estimated residuals based on the finite sequence $\{X_{1:j}^{(n)}\}$ obtained from \eqref{eq:res}.  Assume that $\tilde Z_{j}^{(n)}$ is close to $\hat Z_{j}^{(n)}$ such that for any $\epsilon>0$,
	$$
		\P_n\left(\frac{1}{\sqrt{n}} \sum_{j=1}^n |\hat{Z}_{j}^{(n)} - \tilde{Z}_{j}^{(n)}|^k > \epsilon\right) \to 0, \quad k=1,2.
	$$
\end{enumerate}

\bth~\label{thm:ta}
Assume that \hyperref[cond:n1]{(N1)}, \hyperref[cond:n2]{(N2)}, \hyperref[cond:n3]{(N3)} and \eqref{eq:7} holds, then
$$
	\sup_t\left|\P_n\left(nT_n^{(n)}(h)\le t \right) - \P\left(\|G_h + \xi_h\|_\mu^2\le t\right)\right| \to 0.
$$
\ethe

\bpf[Proof of Theorem~\ref{thm:boot}]
Take $\bbeta^{(n)} = \hat\bbeta$ and $Z_{t}^{(n)} = Z_t^*$.  Here, conditional on the data, $Z_t^*$'s are iid and follow the empirical distribution from $\{\hat{Z}_{1:n}\}$, which converges to the distribution of $Z$ from \eqref{eq:conv:res}. The result follows from Theorem~\ref{thm:ta}.
\epf

\bpf[Proof of Theorem~\ref{thm:ta}]
\noindent
Let $Z_1,Z_2,\ldots$ be a sequence of random variable such that $Z_j\overset{iid}\sim F$.  For each $j$, we have $Z_{j}^{(n)} \cid Z_{j}$. By the Skorohod representation theorem, there exists a sufficiently rich probability space $(\tilde\Omega,\tilde{\mathcal{A}},\tilde{\P})$ where $\tilde\Omega=\{(\omega_1,\omega_2,\ldots):\omega_j\in\Omega_0\}$ for some $\Omega_0$, and functions $z: \Omega_0\to\bbr$, $z^{(n)}: \Omega_0\to\bbr$, such that for each $j$,
$$
	\tilde{Z}_{j}^{(n)} = z^{(n)}(\omega_j) \sim F^{(n)}, \quad \tilde{Z}_{j} = z(\omega_j) \sim F,
$$
and
$$
	\tilde{Z}_{j}^{(n)} \stas \tilde{Z}_j.
$$
This argument is similar to that in \cite{leucht:neumann:2009}.  Since we are only concerned about the distributional limit of $nT_n^{(n)}(h)$, we may assume without loss of generality that $Z_{j}^{(n)}$'s and $Z_j$'s are defined on the same probability space, and $Z_{j}^{(n)} \stas Z_j$ for each $j$. It suffices to prove that in this case,
$$
	nT_n^{(n)}(h) \cid \|G_h + \xi_h\|_\mu^2.
$$

Note that $T_n^{(n)}(h)$ can be written as
$$
	T_n^{(n)}(h) = \int|C_n^{\hat{Z}_n}(s,t)|^2\mu(ds,dt) = \int|C_n^{\hat{Z}_n} - C_n^{Z_n} + C_n^{Z_n}|^2\mu(ds,dt) 
$$
where
$$
	C_n^{\hat{Z}_n}(s,t) := \frac{1}{n} \sum_{j=1}^{n-h} e^{is\hat{Z}_{j}^{(n)}+it\hat{Z}_{j+h}^{(n)}} - \frac{1}{n} \sum_{j=1}^{n-h} e^{is\hat{Z}_{j}^{(n)}} \frac{1}{n} \sum_{j=1}^{n-h} e^{it\hat{Z}_{j+h}^{(n)}}
$$
and
$$
	C_n^{Z_n}(s,t) := \frac{1}{n} \sum_{j=1}^{n-h} e^{isZ_{j}^{(n)}+itZ_{j+h}^{(n)}} - \frac{1}{n} \sum_{j=1}^{n-h} e^{isZ_{j}^{(n)}} \frac{1}{n} \sum_{j=1}^{n-h} e^{itZ_{j+h}^{(n)}}.
$$
The result is proved in two propositions.  In Proposition~\ref{prop:b1}, we show the joint convergence
$$
	(\sqrt{n}C_n^{Z_n}, \sqrt{n}(C_n^{\hat{Z}_n} - C_n^{Z_n})) \cid (G_h,\xi_h),\quad \text{in $\mathcal{C}(K)$},
$$
where $K$ is any compact set in $\bbr^2$. This implies that
$$
	\sqrt{n}C_n^{\hat{Z}_n} \cid G_h + \xi_h,\quad \text{in $\mathcal{C}(K)$}.
$$
Then we justify the convergence of the integral by showing that for any $\varepsilon>0$,
$$
	\lim_{\delta\to0} \limsup_{n\to\infty} \P_n\left( \int_{K_\delta^c}n|C_n^{\hat{Z}_n}|^2 \mu(ds,dt)>\varepsilon \right) = 0.
$$
This is done in Proposition~\ref{prop:b2}.
\epf

\bpr\label{prop:b1}
Given that \hyperref[cond:n1]{(N1)}, \hyperref[cond:n2]{(N2)} and \hyperref[cond:n3]{(N3)} are satisfied, we have
$$
	(\sqrt{n}C_n^{Z_n}, \sqrt{n}(C_n^{\hat{Z}_n} - C_n^{Z_n})) \cid (G_h,\xi_h),\quad \text{in $\mathcal{C}(K)$}.
$$
\epr
\bpf
The proof is divided into the following steps.
\\

\noindent
{\bf Convergence of $C_n^{Z_n}$. }
In this part we show that
$$
	C_n^{Z_n} \cid G_h,\quad \text{in $\mathcal{C}(K)$}.
$$
From Proposition~\ref{prop:joint:conv}, we have $\sqrt{n}C_n^{Z} \cid G_h$, where
$$
	C_n^{Z}(s,t) := \frac{1}{n} \sum_{j=1}^{n-h} e^{isZ_j+itZ_{j+h}} - \frac{1}{n} \sum_{j=1}^{n-h} e^{isZ_{j}} \frac{1}{n} \sum_{j=1}^{n-h} e^{itZ_{j+h}}.
$$ 
It suffices to show that
$$
	\sqrt{n}(C_n^{Z_n} - C_n^{Z}) \cip 0,\quad \text{in $\mathcal{C}(K)$}.
$$
Note that
$$
	C_n^Z(s,t) := \frac{1}{n} \sum_{j=1}^{n-h} U_jV_j - \frac{1}{n} \sum_{j=1}^{n-h} U_j \frac{1}{n} \sum_{j=1}^{n-h} V_j,
$$
where $U_j:=e^{isZ_j}-\varphi_Z(s)$ and $V_j:=e^{itZ_{j+h}}-\varphi_Z(t)$ with $\E U_jV_j = \E U_j = \E V_j=0$.  Similarly,
$$
	C_n^{Z_n}(s,t) := \frac{1}{n} \sum_{j=1}^{n-h} U_{j}^{(n)}V_{j}^{(n)} - \frac{1}{n} \sum_{j=1}^{n-h} U_{j}^{(n)} \frac{1}{n} \sum_{j=1}^{n-h} V_{j}^{(n)},
$$
where $U_{j}^{(n)}(s):=e^{isZ_{j}^{(n)}}-\varphi_{Z^{(n)}}(s)$ and $V_{j}^{(n)}(t):=e^{itZ_{j+h}^{(n)}}-\varphi_{Z^{(n)}}(t)$.  Without loss of generality, here we only show
$$
	\sqrt{n} \left( \frac{1}{n} \sum_{j=1}^{n-h} U_{j}^{(n)} - \frac{1}{n} \sum_{j=1}^{n-h} U_{j}\right) = \frac{1}{\sqrt{n}} \sum_{j=1}^{n-h} (U_{j}^{(n)} - U_j) \cip 0,\quad \text{in $\mathcal{C}(K)$}.
$$
For fixed $s$, the convergence follows since
$$
	\E\left|\frac{1}{\sqrt{n}} \sum_{j=1}^{n-h} (U_{j}^{(n)} - U_j)\right|^2 \le \E|U_{j}^{(n)} - U_j|^2 \to 0,
$$
from bounded convergence.  The finite dimensional convergence can be generalized using the Cram\'er-Wold device.  It remains to prove the tightness of $\frac{1}{\sqrt{n}} \sum_{j=1}^{n-h} (U_{jn} - U_j)$. By equation (7.12) of \cite{billingsley:1999}, the tightness of the process can be implied by
\beao
	\E\left| \frac{1}{\sqrt{n}} \sum_{j=1}^{n-h} (U_{j}^{(n)}(s) - U_j(s)) - \frac{1}{\sqrt{n}} \sum_{j=1}^{n-h} (U_{j}^{(n)}(s') - U_j(s'))\right|^2 \le |s-s'|^{\delta+1} O(1), \quad\text{for some $\delta>0$.}
\eeao
We have
\beao
	&&\E\left| \frac{1}{\sqrt{n}} \sum_{j=1}^{n-h} (U_{j}^{(n)}(s) - U_j(s)) - \frac{1}{\sqrt{n}} \sum_{j=1}^{n-h} (U_{j}^{(n)}(s') - U_j(s'))\right|^2 \\
	&\le& \E\left|U_{j}^{(n)}(s) - U_j(s) - (U_{j}^{(n)}(s') - U_j(s')) \right|^2\\
	&\le& 2\E|e^{isZ_{j}^{(n)}} - e^{is'Z_{j}^{(n)}}|^2 + 2|\varphi_{Z^{(n)}}(s)-\varphi_{Z^{(n)}}(s')|^2+ 2\E|e^{isZ_{j}} - e^{is'Z_{j}}|^2 + 2|\varphi_{Z}(s)-\varphi_{Z}(s')|^2.
\eeao
Note that
$$
	\E|e^{isZ_{j}^{(n)}} - e^{is'Z_{j}^{(n)}}|^2 \le \E|e^{i(s-s')Z_{j}^{(n)}} - 1|^2 \le 2\E|Z_{j}^{(n)}|^2 |s-s'|^2.
$$
The rest of the term can be bounded similarly.  And the tightness is proved.
\\

\noindent
{\bf Convergence of $\sqrt{n}(C_n^{\hat{Z}_n}(s,t) - C_n^{Z_n}(s,t))$. }
In this part we show that
$$
	\sqrt{n}(C_n^{\hat{Z}_n}(s,t) - C_n^{Z_n}(s,t)) \cid \xi_h,\quad \text{in $\mathcal{C}(K)$}.
$$
Similar to \eqref{eq:decomp:en}, we have
\beao
	\sqrt{n}(C_n^{\hat Z_n}(s,t) - C_n^{Z_n}(s,t) )
		= E^{(n)}_n(s,t) - E^{(n)}_n(s,0)\frac{1}{n}\sum_{j=1}^{n-h} \ex^{itZ_{j+h}^{(n)}} - E^{(n)}_n(0,t)\frac{1}{n} \sum_{j=1}^{n-h} \ex^{is\hat Z_{j}^{(n)}}
\eeao
where
$$
	E^{(n)}_n(s,t) :=  \frac{1}{\sqrt{n}} \sum_{j=1}^{n-h} \left( \ex^{is\hat Z_{j}^{(n)}+it\hat Z_{j+h}^{(n)}} -  \ex^{isZ_{j}^{(n)}+itZ_{j+h}^{(n)}} \right).
$$
From the decomposition of $\xi_h$ in \eqref{eq:xi:cal}, it suffices to show that
$$
	E^{(n)}_n(s,t) \cid \bQ^T \mathbf C_h(s,t), \quad \text{in $\mathcal{C}(K)$.}
$$
Uniformly on $(s,t)\in K$, we have
\beao
	E^{(n)}_{n}(s,t)
		&=& \frac{1}{n} \sum_{j=1}^{n-h} \ex^{isZ_j^{(n)}+itZ_{j+h}^{(n)}} (is\sqrt{n}(\hat Z_{j}^{(n)} -  \tilde Z_{j}^{(n)})+it\sqrt{n}(\hat Z_{j+h}^{(n)}- \tilde Z_{j+h}^{(n)})) \\
		&& \quad + \frac{1}{n} \sum_{j=1}^{n-h} \ex^{isZ_{j}^{(n)}+itZ_{j+h}^{(n)}} (is\sqrt{n}(\tilde Z_{j}^{(n)} - Z_{j}^{(n)})+it\sqrt{n}(\tilde Z_{j+h}^{(n)}-Z_{j+h}^{(n)})) + o_p(1) \\
		&=:& E^{(n)}_{n1}(s,t) + E^{(n)}_{n2}(s,t) + o_p(1).
\eeao
From condition \hyperref[cond:n3]{(N3)},
$$
	|E_{n1}^{(n)}(s,t)| \le \frac{|s|+|t|}{\sqrt{n}} \sum_{j=1}^n |\hat Z_{j}^{(n)} -  \tilde Z_{j}^{(n)}| \cip 0, \quad \text{in $\mathcal{C}(K)$.}
$$
It suffices to show that $E^{(n)}_{n2}(s,t)\cid \bQ^T \mathbf{C}_h(s,t)$.  By Taylor expansion,
$$
	E^{(n)}_{n2}(s,t)  =  \sqrt{n}(\hat\bbeta^{(n)} - \bbeta^{(n)})^T \frac{1}{n} \sum_{j=1}^{n-h} \ex^{isZ_{j}^{(n)}+itZ_{j+h}^{(n)}} (is\bl^{(n)}_j(\bbeta^{(n)*}) + it\bl^{(n)}_{j+h}(\bbeta^{(n)*})),
$$
where $\bbeta^{(n)*} = \epsilon\bbeta^{(n)} + (1-\epsilon)\hat\bbeta^{(n)}$ for some $\epsilon\in[0,1]$.
From condition \hyperref[cond:n3]{(N1)}, 
\beqq\label{eq:conv:betaboot}
	\sqrt{n}(\hat\bbeta^{(n)} - \bbeta^{(n)}) \cid \bQ,
\eeqq
follows from the martingale central limit theorem, see Theorem~2 of \cite{scott:1973}. 
It remains to show that
$$
	\frac{1}{n} \sum_{j=1}^{n-h} \ex^{isZ_{j}^{(n)}+itZ_{j+h}^{(n)}} (is\bl^{(n)}_j(\bbeta^{(n)*}) + it\bl^{(n)}_{j+h}(\bbeta^{(n)*})) \cip \mathbf{C}_h(s,t), \quad \text{in $\mathcal{C}(K)$.}
$$
The marginal convergence follows from \hyperref[cond:n2]{(N2)} and the weak law of large number for triangular arrays, see, for example, Theorem~2.2.6 of \cite{durrett:2010}.  The convergence in $\mathcal{C}(K)$ can be extended similar to previous proofs.

\bigskip

\noindent
{\bf Joint convergence of $\sqrt{n}C_n^{Z_n}(s,t)$ and $\sqrt{n}(C_n^{\hat{Z}_n}(s,t) - C_n^{Z_n}(s,t))$. }
The above proofs imply that
$$
	\sqrt{n}C_n^{Z_n}(s,t) - \sqrt{n}C_n^{Z}(s,t) \cip 0, \quad \text{in $\mathcal{C}(K)$,}
$$
and
$$
	\sqrt{n}(C_n^{\hat{Z}_n}(s,t) - C_n^{Z_n}(s,t)) - \sqrt{n}(C_n^{\hat{Z}}(s,t) - C_n^{Z}(s,t)) \cip 0, \quad \text{in $\mathcal{C}(K)$.}
$$
The joint convergence of $\sqrt{n}C_n^{Z_n}(s,t)$ and $\sqrt{n}(C_n^{\hat{Z}_n}(s,t) - C_n^{Z_n}(s,t))$ follows from the joint convergence of $\sqrt{n}C_n^{Z}(s,t)$ and $\sqrt{n}(C_n^{\hat{Z}}(s,t) - C_n^{Z}(s,t))$ in Proposition~\ref{prop:joint:conv}.
\epf

\bpr\label{prop:b2}
For any $\varepsilon>0$,
$$
	\lim_{\delta\to0} \limsup_{n\to\infty} \P\left( \int_{K_\delta^c}n|C_n^{\hat{Z}_n}|^2 \mu(ds,dt)>\varepsilon \right) = 0.
$$
\epr
\bpf
This follows the same steps in the proof of Proposition~\ref{prop:a2} by replacing all $\hat{Z}_j$ with $\hat{Z}_{j}^{(n)}$ and $Z_j$ with $Z_{j}^{(n)}$.
\epf


\section{Proof of Corollary~\ref{thm:arma}} \label{app:arma}

\bpf
In the following we verify conditions \hyperref[cond:m1]{(M1)}, \hyperref[cond:m2]{(M2)}, \hyperref[cond:m3]{(M3)} in Theorem~\ref{thm:meta}.

\noindent
{\hyperref[cond:m1]{(M1)}: }
It can be shown that the pseudo-MLE for $\bbeta$ satisfies the representation in \hyperref[cond:m1]{(M1)}.  We refer to Chapter 10.8 of \cite{brockwell:davis:1991} for details.

\noindent
{\hyperref[cond:m2]{(M2)}: }
From
$$
	Z_t = \frac{\phi(B)}{\theta(B)} X_t =: h(X_{-\infty:t},\bbeta),
$$
we have
$$
	\frac{\partial}{\partial \phi_i} h(X_{-\infty:t},\bbeta) = \frac{B^i}{\theta(B)}X_t = \frac{1}{\theta(B)} X_{t-i}, \quad i=1,\ldots,p,
$$
while
$$
	\frac{\partial}{\partial \theta_i} h(X_{-\infty:t},\bbeta) = \frac{B^j\phi(B)}{(\theta(B))^2}X_t =  \frac{B^j}{\theta(B)}Z_t  = \frac{1}{\theta(B)} Z_{t-j}, \quad j=1,\ldots,q.
$$
Hence
$$
	\bl_0(\bbeta) = \frac{\partial}{\partial\bbeta} h(X_{-\infty:0};\bbeta) = \frac{1}{\theta(B)}(X_{-1},\ldots,X_{-p},Z_{-1},\ldots,Z_{-q})^T.
$$
By the definition of invertibility, there exists a power series for $1/\theta(z)$ such that
$$
	\frac{1}{\theta(z)} = \sum_{j=0}^\infty \xi_j(\bbeta) z^j,
$$
with $\sum_{j=0}^\infty |\xi_j(\bbeta)| < \infty$.
Therefore
$$
	\E\|\bl_0(\bbeta)\|^2 \le p\,\sum_{j=0}^\infty |\xi_j(\bbeta)|^2\E| X_0|^2 + q\,\sum_{k=0}^\infty |\xi_j(\bbeta)|^2\E| Z_0|^2 < \infty.
$$

\noindent
{\hyperref[cond:m3]{(M3)}: }
Note that
$$
	\tilde Z_t - \hat Z_t  = \sum_{j=t}^\infty \pi_j(\hat\bbeta) X_{t-j}.
$$
For $k=1,2$,
\beao
	\frac{1}{\sqrt{n}} \sum_{t=1}^n \left|\tilde Z_t - \hat Z_t \right|^k
	\,\le\, \frac{1}{\sqrt{n}} \sum_{t=1}^n  \sum_{j=t}^\infty \left|\pi_j(\hat\bbeta) X_{t-j} \right|^k\,=\, \sum_{j=0}^\infty |\pi_j(\hat\bbeta) |^k \frac{1}{\sqrt{n}}\sum_{t=1}^{j\wedge n}   \left|X_{t-j}  \right|^k.
\eeao
For any $m<n$,
\beam
	\frac{1}{\sqrt{n}} \sum_{t=1}^n \left|\tilde Z_t - \hat Z_t \right|^k
	\,\le\,\sum_{j=0}^m  |\pi_j(\hat\bbeta) |^k  \frac{1}{\sqrt{n}}\sum_{t=1}^{m} \left|  X_{t-j}  \right|^k + \sum_{j=m+1}^\infty  |\pi_j(\hat\bbeta) |^k  \frac{1}{\sqrt{n}}\sum_{t=1}^{ n} \left| X_{t-j}  \right|^k
	\,=:\,
	I_1+I_2. \label{eq:I1I2}
\eeam
Consider the coefficients $\pi_j(\hat\bbeta)$'s.  By causality, the power series
$$
	\frac{\phi(z)}{\theta(z)} = \sum_{j=0}^\infty \pi_j(\bbeta) z^j
$$
converges for all $|z|<1+\epsilon$ for some $\epsilon>0$.  Then there exists a compact set $\mathbf{C}_\bbeta$ containing $\bbeta$ such that for any $\hat\bbeta \in \mathbf{C}_\bbeta$,
$\sum_{j=0}^\infty \pi_j(\hat\bbeta) z^j$ converges for all $|z|<1+\epsilon/2$.  In particular,
$$
	\pi_j(\hat\bbeta) (1+\epsilon/4)^j \to 0, \quad j\to\infty,
$$
and there exists $K>0$ such that
$$
	|\pi_j(\hat\bbeta)| \le K (1+\epsilon/4)^{-j}.
$$
It follows that for $k=1,2$,
$$
	\sum_{j=0}^\infty  |\pi_j(\hat\bbeta) |^k < \infty
$$
and
$$
	\sum_{j=m}^\infty  |\pi_j(\hat\bbeta) |^k < c(1+\epsilon/4)^{-km}.
$$
Now for \eqref{eq:I1I2}, $I_1$ converges to zero in probability for fixed $m$, while $I_2$ converges to zero uniformly as $m\to\infty$ with order greater than $O(\log(n))$.  This implies that
$$
	\frac{1}{\sqrt{n}} \sum_{t=1}^n \left|\tilde Z_t - \hat Z_t \right|^k \cip 0, \quad k=1,2.
$$

\epf


\section{Proof of Corollary~\ref{thm:garch}} \label{app:garch}

\bpf

In the following we verify conditions \hyperref[cond:m1]{(M1)}, \hyperref[cond:m2]{(M2)}, \hyperref[cond:m3]{(M3)} in Theorem~\ref{thm:meta}.

\noindent
{\hyperref[cond:m1]{(M1)}: }
Given conditions \hyperref[cond:q1]{(Q1)}--\hyperref[cond:q4]{(Q4)}, \cite{berkes:horvath:kokoszka:2003} showed that $\hat\btheta_n$ has limiting distribution
\beao \label{eq:theta:conv}
	\sqrt{n}(\hat\btheta_n - \btheta) 
	\,=\, \frac{1}{\sqrt n} \sum_{t=1}^n \frac{1}{2}(1-Z_t^2)\left\langle\frac{\partial\log\sigma_t^2(\btheta)}{\partial \btheta}, \bB_0^{-1}\right\rangle + o_p(1) 
	\,\cid\, N(\mathbf0,\bB_0^{-1}\bA_0\bB_0^{-1}),
	  \nonumber
\eeao
where
$$
	\bA_0 = \cov \left[ \frac{\partial l_0(\btheta)}{\partial \btheta}\right], \quad \bB_0 = \E \left[ \frac{\partial^2 l_0(\btheta)}{\partial \btheta^2}\right].
$$

\noindent
{\hyperref[cond:m2]{(M2)}: }
We have
$$
	Z_t(\btheta) = h(X_{-\infty:j},\btheta) = \frac{X_t}{\sigma_t(\btheta)},
$$
and
$$
	\bl_0(\btheta) = \frac{\partial}{\partial\btheta} h(X_{-\infty:0};\btheta) = -\frac{X_0}{2\sigma_0^{3}(\btheta)}\frac{\partial\sigma_0^2(\btheta)}{\partial \btheta}
	= -\frac12Z_0\frac{\partial\log\sigma_0^2(\btheta)}{\partial \btheta}.
$$
Lemma~3.1 of \cite{kulperger:yu:2005} showed that
$$
	\E \left( \sup_{\bu\in\Theta} \left|\frac{\partial\log\sigma_t^2(\bu)}{\partial \bu}\right|\right) ^k < \infty,\quad \text{for any $k>0$.}
$$
Hence
\beao
	\E\|\bl_0(\btheta)\|^2 	= \E\left| \frac{1}{2}Z_0 \frac{\partial\log\sigma_0^2(\btheta)}{\partial \btheta} \right|^2 
	\le \frac{1}{4}\left(\E|Z_0|^4   \E\left| \frac{\partial\log\sigma_0^2(\btheta)}{\partial \btheta} \right|^4 \right)^{1/2}
	< \infty.
\eeao

\noindent
{\hyperref[cond:m3]{(M3)}: }
Theorem~1.3 and Lemma~3.5 of \cite{kulperger:yu:2005} show, respectively, that
$$
	\frac{1}{\sqrt{n}} \sum_{t=1}^{n} |\hat Z_t - \tilde Z_t| = o_p(1),
$$
and
$$
	\sum_{t=1}^{n} |\hat Z_t - \tilde Z_t| = O_p(1).
$$
Hence 
$$
	\frac{1}{\sqrt{n}} \sum_{t=1}^{n} |\hat Z_t - \tilde Z_t|^2 \le \frac{1}{\sqrt{n}} \sum_{t=1}^{n} |\hat Z_t - \tilde Z_t|\, \sum_{t=1}^{n} |\hat Z_t - \tilde Z_t| = o_p(1).
$$

\epf

\end{document}